\documentclass[journal]{IEEEtran}
\usepackage{amsfonts}
\usepackage{bm}
\usepackage{amsmath,amsfonts,amsthm,amssymb}
\usepackage{setspace}
\usepackage{Tabbing}
\usepackage{fancyhdr}
\usepackage{lastpage}
\usepackage{extramarks}
\usepackage{chngpage}
\usepackage{soul,color}
\usepackage{epsfig}
\usepackage{graphicx,float,wrapfig,subfigure}
\usepackage{longtable}
\usepackage{wrapfig}
\usepackage{stfloats}
\usepackage{cite}
\usepackage{graphicx}
\usepackage{array}
\usepackage{multirow}{\tiny }
\usepackage{multicol}
\usepackage{colortbl}
\usepackage{tabularx}
\usepackage{mdwmath}
\usepackage{mdwtab}
\usepackage{color}
\usepackage{verbatim}
\usepackage{footnote}
\makesavenoteenv{table}
\usepackage{amsmath}
\usepackage{tikz}
\usetikzlibrary{calc}
\usepackage{flowchart}
\usetikzlibrary{shapes.geometric, arrows}
\usepackage{overpic}
\usepackage{url}
\usepackage[colorlinks,linkcolor=black,anchorcolor=black,citecolor=black,urlcolor=black]{hyperref} 
\usepackage{breakurl}
\usepackage[flushleft]{threeparttable}
\usepackage{comment}
\allowdisplaybreaks

\usepackage{xr}
\makeatletter
\newcommand*{\addFileDependency}[1]{% argument=file name and extension
  \typeout{(#1)}% latexmk will find this if $recorder=0 (however, in that case, it will ignore #1 if it is a .aux or .pdf file etc and it exists! if it doesn't exist, it will appear in the list of dependents regardless)
  \@addtofilelist{#1}% if you want it to appear in \listfiles, not really necessary and latexmk doesn't use this
  \IfFileExists{#1}{}{\typeout{No file #1.}}% latexmk will find this message if #1 doesn't exist (yet)
}
\makeatother

\tikzstyle{startstop} = [rectangle, rounded corners, minimum width=1.5cm, minimum height=0.5cm,text centered, draw=black, fill=red!30]
\tikzstyle{io} = [trapezium, trapezium left angle=70, trapezium right angle=110, minimum width=1.5cm, minimum height=0.5cm, text centered, draw=black, fill=blue!30]
\tikzstyle{process} = [rectangle, minimum width=1.5cm, minimum height=0.5cm, text centered, draw=black, fill=orange!30]
\tikzstyle{decision} = [diamond,  minimum width=2.5cm, minimum height=0.5cm, text centered, draw=black, fill=green!30]
\tikzstyle{arrow} = [thick,->,>=stealth]
\makeatletter
\def\env@cases{%
  \let\@ifnextchar\new@ifnextchar
  \left\{
  \def\arraystretch{1.2}%
  \array{@{}l@{\,\,}l@{}}%
}%
\makeatother

\includecomment{extended}
\excludecomment{reduced}

\newcounter{remark}
\newenvironment{remark}[1][]{\refstepcounter{remark}\par
	\textbf{Remark~\theremark #1:} \rmfamily}

\newcounter{proposition}
\newenvironment{proposition}[1][]{\refstepcounter{proposition}\par
	\textbf{Proposition~\theproposition #1:} \rmfamily}

\newcounter{definition}
\newcounter{theorem}
\newenvironment{theorem}[1][]{\refstepcounter{theorem}\par
	\textbf{Theorem~\thetheorem #1:} \rmfamily}

\newcounter{assumption}
\newcounter{corollary}
\newenvironment{corollary}[1][]{\refstepcounter{corollary}\par
	\textbf{Corollary~\thecorollary #1:} \rmfamily}

\newcounter{lemma}
\newenvironment{lemma}[1][]{\refstepcounter{lemma}\par
	\textbf{Lemma~\thelemma #1:} \rmfamily}

\newcounter{algorithm}
\graphicspath{{figures/}}

\begin{document}

\title{
    Optimal Control of AGC Systems Considering Non-Gaussian Wind Power Uncertainty
	}

\author{
Xiaoshuang~Chen,~\IEEEmembership{Student Member,~IEEE,}
Jin~Lin,~\IEEEmembership{Member,~IEEE,}
Feng~Liu,~\IEEEmembership{Member,~IEEE,}
and~Yonghua~Song,~\IEEEmembership{Fellow,~IEEE}

\thanks{
%	This work was supported by National High-Technology Research and Development Program (“863” Program) of China (2014AA051901), National Natural Science Foundation of China (51577096,51761135015), International S\&T Cooperation (2016YFE0102600).
	X. Chen, J. Lin, F. Liu and Y. Song are with Department of Electrical Engineering, Tsinghua University, Beijing 100084, China(email: linjin@tsinghua.edu.cn).
%    X. Chen, J. Lin and F. Liu are with the State Key Laboratory of Control and Simulation of Power Systems and Generation Equipment, Department of Electrical Engineering, Tsinghua University, Beijing 100084, China.% (email: linjin@tsinghua.edu.cn).

%    Y. Song is with the Department of Electrical and Computer Engineering, University of Macau, Macau, China, and the Department of Electrical Engineering, Tsinghua University, Beijing 100084, China.
	}
}

\maketitle
\begin{abstract} 
    Wind power uncertainty poses significant challenges for automatic generation control (AGC) systems. It can enhance control performances to explicitly consider wind power uncertainty distributions within controller design. However, widely accepted wind uncertainties usually follow non-Gaussian distributions which may lead to complicated stochastic AGC modeling and high computational burdens. To overcome the issue, this paper presents a novel It\^{o}-theory-based model for the stochastic control problem (SCP) of AGC systems, which reduces the computational burden of optimization considering non-Gaussian wind power uncertainty to the same scale as that for deterministic control problems. We present an It\^{o} process model to exactly describe non-Gaussian wind power uncertainty, and then propose an SCP based on the concept of stochastic assessment functions (SAFs). Based on a convergent series expansion of the SAF, the SCP is reformulated as a certain deterministic control problem without sacrificing performance under non-Gaussian wind power uncertainty. The reformulated control problem is proven as a convex optimization which can be solved efficiently. A case study demonstrates the efficiency and accuracy of the proposed approach compared with several conventional approaches.
\end{abstract}
%\vspace{-3mm}

\vspace{-2mm}
\section{Introduction} \label{section:introduction}
\IEEEPARstart{I}{n} past few years, the penetration of wind generations has been increasing rapidly \cite{Markorov2010Analysis}. However, the wind generations are introducing much more uncertainty into power systems, and this uncertainty has severe impacts on the automatic generation control (AGC) systems \cite{Wang2016Fully,Wan2014Optimal}. The uncertainty brought by wind generations is usually non-Gaussian, which is challenging for uncertainty modeling and AGC design under this non-Gaussian uncertainty.

%{\color{red}two approaches: AGC and scenario-based programming. Add AGC papers in the introduction. discuss the problem in AGC rather than control.}

Techniques for handling wind power uncertainty in AGC systems include proportional-integral (PI) control \cite{Sahu2014Hybrid,Xu2016Dynamic}, robust control \cite{Yao2017Robust, Sharma2016Robust}, and model predictive control (MPC) \cite{Ma2017Distributed, Ersdal2016Model, Cominesi2018Two, Mcnamara2018Model, Jiang2016Explicit}. PI AGC controllers are among the most popular controllers in practice, and there have been studies on several variants thereof, such as fuzzy PI controllers \cite{Sahu2014Hybrid} and dynamic-gain PI controllers \cite{Xu2016Dynamic}. In robust AGC approaches, the uncertainty is modeled using intervals or bounded sets, and the controllers are designed to guarantee the performance in the worst case scenario\cite{Yao2017Robust, Sharma2016Robust}.
%MPC approaches outperform the abovementioned approaches when in the optimal control problem in AGC systems with state constraints. Unlike the close-loop implementations of PI controllers and robust controllers, 
In MPC AGC controllers\cite{Ma2017Distributed, Ersdal2016Model, Cominesi2018Two, Mcnamara2018Model, Jiang2016Explicit}, an open-loop  optimal control problem is solved repeatedly in a receding-horizon manner.
%It does not model the distribution of the uncertainty. 
In this approach, the uncertainty is modeled as a deterministic input based on the calculated prediction, and the prediction is revised over the control period. Therefore, MPC requires only the solution of a deterministic optimization problem at each time step.
Although the approaches mentioned above achieve success for many AGC problems, a key issue is that the probability distribution of the uncertainty is not explicitly considered in the controller model, which may limit the control performance. For example, MPC uses a receding-horizon implementation to cope with the uncertainty. To achieve better performance, a larger number of prediction steps and a smaller time step are needed, which will result in a heavy online computational burden \cite{Camacho2013Model}.

Problems that explicitly consider the probability distribution of the uncertainty in control systems are usually regarded as stochastic control problems (SCPs).
%Generally, The optimal SCP can be solved via stochastic dynamic programming (SDP)\cite{Oksendal2003Stochastic}. However, SDP suffers from the curse of dimensionality and is practical only for very small systems (such as 2nd- or 3rd-order systems \cite{Picarelli2015Stochastic}).
Generally, optimal SCPs with Gaussian distributions can be solved via quadratic programming \cite{Hadjiyiannis2011Efficient}. However, the wind power uncertainty in AGC systems is usually non-Gaussian. In previous research, the long-tailed characteristic of the wind power uncertainty has usually been modeled with a beta distribution or Laplace distribution \cite{Bludszuweit2008Statistical}. SCPs with non-Gaussian uncertainty are usually solved via scenario-based stochastic programming (SBSP) algorithms, in which the uncertainties are characterized by a finite set of random realizations and the SCP is solved via a deterministic optimization problem \cite{Yang2016Stochastic}. The SBSP approach is widely used in intra-day applications of power systems\cite{Fu2016Multiobjective, Yang2016Stochastic}. Although a few studies have applied SBSP for the assessment and control of AGC systems\cite{Engels2017Combined, Donadee2012Stochastic,Gross2001Analysis,Chang2014Modeling}, it is still challenging to do so because the number of variables is proportional to the number of scenarios, which is usually large to ensure good control performance \cite{Campi2008Exact}. Therefore, this approach is not typically applicable in the AGC context.

In summary, several challenges arise when dealing with uncertainty in AGC systems: 1) The uncertainty is usually non-Gaussian and dependent on time\cite{Tabone2015Modeling,Bludszuweit2008Statistical}.
%Algorithms for coping with non-Gaussian uncertainty usually have a high order of complexity compared with the algorithms for deterministic systems\cite{Mesbah2016Stochastic,Picarelli2015Stochastic}.
2) It is challenging to explicitly consider uncertainty in AGC systems while incurring a feasible computational burden. 
Against this background, this paper proposes a novel It\^{o}-theory-based approach, of which the key idea is to introduce the concept of stochastic assessment functions (SAFs) to express the expectation values in SCPs and then transform the SAFs into a few deterministic assessment functions (DAFs) to significantly reduce the computational burden.

The major contributions of this paper include the following:
\begin{enumerate}
\item An It\^{o} model of an AGC system considering the the non-Gaussian stochastic characteristics of wind farms is presented, in which the wind power uncertainty is modeled in the form of It\^{o} processes and the AGC system is described by a stochastic differential equation (SDE) model.
%%%Editor - Acronyms and abbreviations are often defined at their first
%%%appearance in the text and then used throughout the remainder of the
%%%manuscript. Please consider adhering to this convention.
This formulation is applicable for various distributions of wind power uncertainty.
	\item Based on the It\^{o}-AGC system model, an SCP is formulated in which the expectation values are expressed as SAFs. We then reveal the SAFs into a series of DAFs with the proof of convergence. Unlike in the SBSP approach, no scenario is needed in this approach.
	\item A convex optimization formulation of the It\^{o}-AGC is formulated based on the series expansion of the SAFs, which consequently achieves a good trade-off between the accuracy and the computational burden.
\end{enumerate}

Following this introduction, Section \ref{section:model} presents the modeling of wind power uncertainty and AGC systems. Section \ref{section:scp} formulates the SCP for an AGC system via the concept of SAFs. Section \ref{section:convex} presents the series expansion of an SAF and the transformation of the SCP into a convex optimization problem. A case study is presented in Section \ref{section:case}, and Section \ref{section:conclusion} concludes the paper. 
%\begin{reduced}
%Due to space constraints, some proofs are put in the appendix which is accessible online \cite{Chen2018Stochastic}.
%\end{reduced}

\color{black}
\vspace{-2mm}
\section{It\^{o} Modeling of an AGC System} \label{section:model}
This section presents an It\^{o} model of an AGC system. An It\^{o} process model is used to describe the stochastic characteristics of the sources of uncertainty, based on which the dynamics of the AGC system are modeled as an SDE. All variables are expressed as deviations from the operating point.

\vspace{-2mm}
\subsection{Power System Model} \label{section:model-agc}
\vspace{-1mm}
\subsubsection{Generators}
A typical synchronous generator consists of a governor, a turbine and a rotating mass. Let $\Omega^G$ be the set of generators; then, for $i\in\Omega^G$, we have 
\vspace{-2mm}
\begin{equation} \label{eq:model-generator-dynamic}
\begin{split}
\dot{P}_{i}^g &= -\frac{1}{T_{i}^g}\left(P_{i}^g + \frac{\omega_i}{R_i}-P_{i}^{ref}\right) \\
\dot{P}_{i}^m &= -\frac{1}{T_i^t}\left(P_i^m-P_i^g\right) 
%\vspace{-1mm}
\\
%\vspace{-1mm}
\dot{w}_i &= -\frac{1}{H_i}\left(D_i\omega_i - P_i^m + Z^{(i)} + \sum_{l\in Adj(i)}P_{l}^L\right)
\end{split}
\end{equation}
where $P_i^g$ and $P_i^m$ are the outputs of the governor and turbine, respectively; $T_i^g$ and $T_i^t$ are the time constants of the governor and turbine, respectively; $R_i$, $H_i$, $D_i$ and $\omega_i$ are the droop, the inertia, the damping, and the speed of the rotating mass, respectively; $Z^{(i)}$ is the wind power uncertainty at bus $i$, which will be discussed in Section \ref{section:uncertainty}; $P_i^{ref}$ is the power reference, which is controllable ; $P_{ij}$ is the power flow from bus $i$ to bus $j$;
%%%Editor - Please ensure that the intended meaning has been maintained in the
%%%above edit. Alternatively, you may prefer ``the power flow between bus $i$
%%%and bus $j$''.
and $Adj(i)$ is the set of buses adjacent to bus $i$.

In secondary frequency control, the time step is usually several seconds. Therefore, it is usually assumed that the time constant of primary frequency control is negligible, i.e., $T_i^g=T_i^t=0$ \cite{Jiang2016Explicit}; consequently, we have
\vspace{-2mm}
\begin{equation}
\begin{split}
P_{i}^g &= P_i^m = P_{i}^{ref}-\omega_i/R_i \\
\end{split}
\end{equation}
Therefore, \eqref{eq:model-generator-dynamic} can be simplified as
\vspace{-2mm}
\begin{equation} \label{eq:model-generator-dynamic-simplified}
\begin{split}
\dot{\omega}_i &= -\frac{1}{H_i}\left[\left(D_i+\frac{1}{R_i}\right)\omega_i - P_i^{ref} - Z^{(i)} + \sum_{j\in Adj(i)}P_{ij}\right]
\end{split}
\end{equation}

\subsubsection{Branch Flow}
The branch flow $P_{ij}$ satisfy
\vspace{-1mm}
\begin{equation} \label{eq:model-network-dcflow}
\dot{P}_{ij} = B_{ij}(\omega_i-\omega_j)
\vspace{-1mm}
\end{equation}
with $B_{ij}$ being a constant defined as
\vspace{-1mm}
\begin{equation}
B_{ij} = \frac{|V_i||V_j|}{x_{ij}}\cos\left(\theta_i^0-\theta_j^0\right)
\vspace{-1mm}
\end{equation}
where $V_i$ and $V_j$ are the nominal bus voltages, $x_{ij}$ is the line reactance, and $\theta_i^0$ and $\theta_j^0$ are the operating points of the phase angles. The same model can be found in \cite{Li2016Connecting}. Note that \eqref{eq:model-network-dcflow} omits the initial deviation $P_{ij}(0)$ in the branch flow. In practice, $P_{ij}(0)$ cannot be arbitrary but instead must satisfy
\vspace{-1mm}
\begin{equation}
P_{ij}(0) = B_{ij}\left(\theta_i(0)-\theta_j(0)\right)
\vspace{-1mm}
\end{equation}
for some vector $\bm{\theta}$ \cite{Zhao2014Design}.
%Though \eqref{eq:model-network-dcflow} is based on the DC power flow model, it should be emphasized that other linear power flow models, such as the models in \cite{Wan2017Maximum,Xing2017Model}, can also be applied in the proposed model.

%\subsubsection{Storage}
%Recently, storage systems are widely considered in AGC systems \cite{Chakraborty2018Automatic}. Denote $\Omega^{E}$ as the set of storage systems. For $i\in \Omega^E$, the dynamics of the $i$-th storage system can be described as\cite{Qin2016Online}:
%%\vspace{-0.5mm}
%\begin{equation} \label{eq:model-eu}
%\dot{E}_i^E = -\alpha_i^E E_i^E - P_i^E
%\end{equation}
%%\vspace{-0.5mm}
%where $E_i^E$ is the stored energy, $P_i^E$ is the output power, and $\alpha_i^E$ is the dissipation factor.

\subsubsection{AGC Signals}
In an AGC system \cite{Gross2001Analysis}, the frequency deviation $\Delta_f$ is the weighted average of the local frequencies:
\vspace{-2mm}
\begin{equation} \label{eq:model-frequency-deviation}
\Delta_f = \frac{1}{2\pi}\frac{\sum_{i\in\Omega^G}H_i\omega_i}{\sum_{i\in\Omega^G}H_i}
\end{equation}

\vspace{-1mm}
The area control error (ACE) of area $m$ is defined as
\vspace{-1mm}
\begin{equation} \label{eq:model-ace}
ACE_m = \sum_{ij\in \Omega^L_m}P_{ij}-b_m\Delta_f
\vspace{-1mm}
\end{equation}
where $ACE_m$ is the ACE of area $m$, $\Omega^L_m$ is the set of tie-lines connected to area $m$, and $b_m$ is the bias factor of area $m$.

\vspace{-2mm}
\subsection{Wind Uncertainty Modeling Based on It\^{o} Process} \label{section:uncertainty}
The distribution of the wind power uncertainty in an AGC system may be non-Gaussian in practice, and thus, the model of the uncertainty sources should support different kinds of probability distributions.

We use an It\^{o} process model to express the stochastic characteristics of the wind power uncertainty. Specifically, let $Z_t$ be the wind power uncertainty (as denoted in \eqref{eq:model-generator-dynamic}, but the bus number $i$ is omited for convenience); then, its characteristics are described by the following SDE:
\vspace{-1mm}
\begin{equation} \label{eq:sc-source}
dZ_t = \mu(Z_t)dt + \sigma(Z_t)dW_t
\end{equation}
where $\mu$ is the drift function that describes the deterministic characteristics, whereas $\sigma$ is the diffusion function that describes the stochastic characteristics. $W_t$ is the integrand, which is a standard Wiener process\cite{Pardoux2014Stochastic}. 

%In some literatures \cite{Apostolopoulou2016Assessment, Yuan2015Stochastic}, $\mu(\cdot)$ is set affine to $Z_t$ and $\sigma(\cdot)$ is a constant. In such case, $Z_t$ satisfies a Gaussian distribution. An arbitrary probability distribution function (PDF) can be described \eqref{eq:sc-source} with suitable $\mu$ and $\sigma$. 

With different specifications for $\mu(\cdot)$ and $\sigma(\cdot)$, the It\^{o} process model can describe different probability density functions (PDFs). Here, we present the $\mu(\cdot)$ and $\sigma(\cdot)$ specifications for various distributions of the wind power uncertainty, as shown in Table \ref{tab:ito-distribution}.
\begin{extended}
	The proof is given in Appendix \ref{appendix:ito-distribution}, which also shows that an It\^{o} process model can describe the uncertainties of arbitrary PDFs.
\end{extended}
\begin{reduced}
	The proof is given in Appendix A of the extended version \cite{Chen2018Stochastic}, which also shows that an It\^{o} process can describe uncertainties of arbitrary PDFs.
\end{reduced}

The It\^{o} process model provides an SDE form for an AGC system, which is important in the derivations throughout the remainder of this paper. For this purpose, we use an It\^{o} process for AGC modeling, as introduced in the next subsection.
\begin{table}[!t]
	\renewcommand{\arraystretch}{1.1}
	%	\linespread{1.5}
	\centering
	\begin{small}
		\caption{Specific Formulations for Some Typical Distributions}
		\vspace{-2mm}
		\begin{tabular}{cccc}
			\hline\hline
			Distribution&PDF&$\mu(z)$&$\sigma^2(z)$\\
			\hline
			Gaussian&$N(a,b)$&$-(z-a)$&$2b$\\
			\hline
			Beta&$B(a,b)$&$-(z-\frac{a}{a+b})$&$\frac{2}{a+b}z(1-z)$\\
			\hline
			Gamma&$\Gamma (a, b)$&$-(z-a/b)$&$2z/b$\\
			\hline
			Laplace&$L(a,b)$&$-(z-a)$&$2b|z-a|+2b^2$ \\
			%			\hline
			%			Weibull &$W(\lambda, k)$&$-(z-\lambda\Gamma(1+\frac{1}{k}))$&$\frac{2}{p(z)}\left[\lambda\Gamma(1+\frac{1}{k},\left(\frac{z}{\lambda}\right)^k)-\lambda\Gamma(1+\frac{1}{k})-e^{-(z/\lambda)^k}\right]$ \\
			\hline \hline
		\end{tabular}\label{tab:ito-distribution}
		\vspace{-5mm}
	\end{small}
\end{table}
%
%According to the steady-state Fokker-Planck equation \cite{Pardoux2014Stochastic}, we have
%\begin{equation} \label{eq:fokker-planck}
%\frac{1}{2}\frac{\partial^2 \left(\sigma^2(z)p(z)\right)}{\partial z^2} = \frac{\partial \left(\mu(z)p(z)\right)}{\partial z}
%\end{equation}
%By integrating \eqref{eq:fokker-planck} over $z$, we have
%\begin{equation} \label{eq:relationship-mu-sigma}
%\frac{1}{2}\sigma^2(z)p(z) = \int_{\infty}^z\mu(z')p(z')dz'
%\end{equation}
%Where the boundary condition at infinity is used. In practice, the PDF $p$ can be obtained via statistical analysis of historical data. Once $p$ is obtained, the relationship is described by \eqref{eq:relationship-mu-sigma}. Therefore, the It\^{o} process enables us to model uncertainty sources with arbitrary PDFs.
%
%Moreover, since $\mu$ and $\sigma$ are not uniquely determined by \eqref{eq:relationship-mu-sigma}, one can fix $\mu$ as a predefined function, such as an affine function with respect to $z$, and then compute $\sigma$ according to \eqref{eq:relationship-mu-sigma}.

\color{black}
\vspace{-2mm}
\subsection{AGC Modeling Based on It\^{o} Process}
\subsubsection{Variable Sets}
Now, we transform the models presented above into an SDE form. We first define some variable sets.
\begin{itemize}
	\item Uncertain wind generations, denoted by $\bm{Z}_t$. Each element of the vector $\bm{Z}_t$ satisfies \eqref{eq:sc-source}. For simplicity, we use the following matrix form:
	\vspace{-1mm}
	\begin{equation} \label{eq:model-z}
	d\bm{Z}_t = \mu(\bm{Z}_t)dt + \sigma(\bm{Z}_t)d\bm{W}_t
	\vspace{-1mm}
	\end{equation}
	\item Variables satisfying ordinary differential equations (ODEs), including $\omega_i$ for $i\in \Omega^G$ and $P_{ij}$. These variables are denoted by the vector $\bm{X}_t$.
	\item Variables satisfying algebraic equations, including $ACE_m$ for area $m$ and $\Delta_f$. These variables are denoted by the vector $\bm{Y}_t$. Note that $\bm{Y}_t$ does not need to exist in the final model since it can be substituted with $\bm{X}_t$ and $\bm{Z}_t$ in accordance with \eqref{eq:model-frequency-deviation} and \eqref{eq:model-ace}.
	\item Control variables, i.e., $P_{i}^{ref}$ for $i \in \Omega^G$. These variables are denoted by $\bm{U}_t$.
\end{itemize}
Thus, the ODEs in \eqref{eq:model-generator-dynamic-simplified} \eqref{eq:model-network-dcflow} can be rewritten as:
\vspace{-1mm}
\begin{equation} \label{eq:model-x}
\dot{\bm{X}}_t = \bm{AX}_t + \bm{BU}_t + \bm{CZ}_t
\vspace{-1mm}
\end{equation}
where $\bm{A}$, $\bm{B}$ and $\bm{C}$ are matrix coefficients. This linear-differential-equation formulation has been widely used in previous studies\cite{Ma2017Distributed, Ersdal2016Model, Jiang2016Explicit}.

\subsubsection{Control Policy}
Generally, there are two kinds of control policies: one is state feedback control \cite{Cannon2009Probabilistic,Skaf2010Design}, i.e.,
\vspace{-1mm}
\begin{equation} \label{eq:state-control}
\bm{U}_t = u(\bm{X}_t, \bm{Z}_t)
\vspace{-1mm}
\end{equation}
and the other is disturbance feedback control \cite{Chatterjee2011Stochastic,Ding2016Optimal,Oldewurtel2008Tractable,Korda2011Strongly}, i.e.,
\vspace{-1mm}
\begin{equation} \label{eq:disturbance-control}
\bm{U}_t = u(\bm{Z}_t)
\vspace{-1mm}
\end{equation}
Since the performances of these two types of control policies have been proven to be equivalent \cite{Skaf2010Design}, we adopt the disturbance feedback control policy \eqref{eq:disturbance-control} since it will be helpful for establishing the convex optimization problem. Moreover, we denote the set of control policies by $\mathcal{U}$.

\subsubsection{The It\^{o}-AGC Model}
%According to Section \ref{section:model-agc}, $\bm{X}_t$, $\bm{Y}_t$, $\bm{Z}_t$ and $\bm{U}_t$ satisfy the following equations
%\begin{equation} \label{eq:ito-whole}
%\begin{split}
%\bm{\dot{X}}_t &= \bm{A}_X\bm{X}_t + \bm{A}_Y\bm{Y}_t + \bm{A}_Z\bm{Z}_t + \bm{A}_U\bm{U}_t \\
%\bm{Y}_t &= \bm{B}_X\bm{X}_t + \bm{B}_Z\bm{Z_t} + \bm{B}_U\bm{U}_t \\
%d\bm{Z}_t &= \mu(\bm{Z}_t)dt + \sigma(\bm{Z}_t)d\bm{W}_t
%\end{split}
%\end{equation}

According to \eqref{eq:model-z} and \eqref{eq:model-x}, the \textit{stochastic system} (SS) of interest can be formulated as

\noindent\textbf{SS}:
\vspace{-2mm}
\begin{equation} \label{eq:ito-reduced}
\begin{split}
d\bm{X}_t &= (\bm{A}\bm{X}_t + \bm{B}u(\bm{Z}_t) + \bm{C}\bm{Z}_t)dt\\
d\bm{Z}_t &= \mu(\bm{Z}_t) + \sigma(\bm{Z}_t)d\bm{W}_t \\
\bm{X}_0&=\bm{x}_0,\bm{Z}_0=\bm{z}_0
\end{split}
\end{equation}
where $\bm{x}_0$ and $\bm{z}_0$ are the initial values of $\bm{X}_t$ and $\bm{Z}_t$, respectively. For convenience, let $N_x$ be the dimensionality of $\bm{X}_t$, and let $N_z$ be the dimensionality of $\bm{Z}_t$. The major difference between \eqref{eq:ito-reduced} and existing ODE-based AGC models \cite{Ma2017Distributed, Ersdal2016Model, Jiang2016Explicit} is the model of the disturbance $\bm{Z}_t$, which is an It\^{o} process model here. Therefore, the model used here is an SDE model rather than an ODE model. This SDE model describes both the stochastic characteristics and the dynamics of the SS in a unified way; moreover, it enables an analytical formulation of the SCP
%%%Editor - Once an acronym or abbreviation has been defined, it is often used
%%%throughout the remainder of the manuscript. Please consider adhering to this
%%%convention.
without scenario generation, as will be discussed in the remainder of this paper.

\vspace{-2mm}
\section{Stochastic Optimal Control of It\^{o}-AGC} \label{section:scp}
This section presents the SAFs to describe the objective function and constraints based on the model of It\^{o}-AGC and then formulates the SCP that will be studied in this paper.

\subsection{Objective}
The objective function for AGC \cite{Jiang2016Explicit} is
\vspace{-2mm}
\begin{equation} \label{eq:objective}
\begin{split}
J =& \mathbb{E}^{\bm{x}_0,\bm{z}_0}\left[\int_0^T \left(\lambda_{ACE}ACE_m^2(t) +\bm{U}_t^{\top}\bm{\Lambda}_U\bm{U}_t\right)dt\right] \\
&+\mathbb{E}^{\bm{x}_0,\bm{z}_0}\left[\mu_{ACE}ACE_m^2(T)\right]
\end{split}
\end{equation}
where $T$ is the terminal time, $\lambda_{ACE}$ is the coefficient of the ACE, $\Lambda_U$ is the coefficient of the control variable, and $\mu_{ACE}$ is the coefficient of the terminal state. $\mathbb{E}^{\bm{x}_0,\bm{z}_0}\{ \cdot \}$ denotes the expectation operator under the initial conditions $\bm{x}_0$ and $\bm{z}_0$, and the superscript ``${\top}$'' denotes the transpose operation. The objective function expressed in \eqref{eq:objective} shows the trade-off between the control performance, as measured by $ACE_m$, and the cost, as measured by $\bm{U}_t$.

For convenience, we use the following general form to represent the objective function:
\vspace{-1mm}
\begin{equation} \label{eq:objective-general}
\min J = \mathbb{E}^{\bm{x}_0,\bm{z}_0}\left[\int_0^Tf(\bm{S}_t)dt + g(\bm{S}_T)\right]
\vspace{-1mm}
\end{equation}
where $\bm{S}_t$ represents the concatenation of $\bm{X}_t$, $\bm{U}_t$, and $\bm{Z}_t$:
\begin{equation} \label{eq:s}
\bm{S}_t = 
\begin{bmatrix}
\bm{X}_t \\
u(\bm{Z}_t) \\
\bm{Z}_t
\end{bmatrix}
\end{equation}

It is clear that \eqref{eq:objective-general} is a generalization of \eqref{eq:objective} since the ACE can be expressed as a linear combination of the state variables $\bm{X}_t$, $\bm{U}_t$ and $\bm{Z}_t$. Moreover, $f$ and $g$ are both convex and of up to quadratic order.

\vspace{-2mm}
\subsection{Constraints of an AGC System}
The constraints in an AGC system \cite{Jiang2016Explicit} include:
\vspace{-1mm}
\begin{equation} \label{eq:model-constraint}
%\begin{array}{c}
\begin{split}
&-\bar{P}_i^G \leq P_i^G \leq \bar{P}_i^G, \forall i \in \Omega^G \\
%-\bar{P}_i^E < P_i^E \leq \bar{P}_i^E, - \bar{E}_i^E \leq E_i^E \leq \bar{E}_i^E, \forall i \in \Omega^E \\
&-\bar{P}_{ij} \leq P_{ij} \leq \bar{P}_{ij}, \forall i,j \in \Omega^G, i \in Adj(j) \\
&-\bar{\Delta}_f \leq \Delta_f \leq \bar{\Delta}_f
\end{split}
\end{equation}

All of these constraints are linear constraints, and according to the models presented in Section \ref{section:model-agc}, the variables in these constraints can be regarded as linear combinations of $\bm{S}_t$; therefore, these constraints can be uniformly rewritten as:
\vspace{-1mm}
\begin{equation} \label{eq:model-constraint-uniform}
\bm{\phi}_i^{\top} \bm{S}_t \leq \bar{\phi}_i, \forall i \in \Omega^C
\vspace{-1mm}
\end{equation}
where $\Omega^C = \Omega^G \cup \Omega^E \cup \Omega^L$, 
$\bm{\phi}_i$ is the coefficient vector specifying the linear combination of $\bm{S}_t$, and $\bar{\phi}_i$ is the upper bound on the $i$-th constraint.

The almost-sure constraint expressed in \eqref{eq:model-constraint-uniform} can be naturally relaxed to a chance constraint \cite{Calafiore2006Distributionally}:
\vspace{-1mm}
\begin{equation} \label{eq:chance-constraint}
\Pr \left\{\bm{\phi}_i^{\top}\bm{S}_t \leq \bar{\phi}_i\right\} \geq \gamma
\vspace{-1mm}
\end{equation}

According to \cite{Calafiore2006Distributionally}, the chance constraint given in \eqref{eq:chance-constraint} has a tractable inner approximation:
\begin{equation} 
\kappa_\gamma\sqrt{\text{var}\{\bm{\phi}_i^{\top}\bm{S}_t\}} + \mathbb{E}^{\bm{x}_0,\bm{z}_0}\left\{\bm{\phi}_i^{\top}\bm{S}_t\right\}\leq \bar{\phi}_i \label{eq:sc-constraint-tractable}
\end{equation}
where
\vspace{-1mm}
\begin{equation} \label{eq:sc-constraint-var}
\text{var}\{\bm{\phi}_i^{\top}\bm{S}_t\} = \mathbb{E}^{\bm{x}_0,\bm{z}_0}\bm{\phi}_i^{\top}\bm{S}_t\bm{S}_t^{\top}\bm{\phi}_i-\left(\mathbb{E}^{\bm{x}_0,\bm{z}_0}\left\{\bm{\phi}_i^{\top}\bm{S}_t\right\}\right)^2
\vspace{-1mm}
\end{equation}

\subsection{Stochastic Control Problem}
Before presenting the SCP, it is necessary to introduce an SAF to describe the objective function and the constraints in a unified form. Specifically, given that $\alpha$ and $\beta$ are two functions of $\bm{S}_t$, we define the SAF formulation for the SS in \eqref{eq:ito-reduced} as:
\begin{equation} \label{eq:sap}
\begin{split}
\mathcal{P}_{\alpha, \beta}(t,\bm{x}_0, \bm{z}_0)= \mathbb{E}^{\bm{x}_0, \bm{z}_0}\left\{{\int_0^t\alpha(\bm{S}_s)ds} + \beta(\bm{S}_t)\right\}
\end{split}
\end{equation}
where the subscript $s$ is used to replace the time variable $t$, which has been used as the upper endpoint of the integral.

Via the SAF concept, it is possible to express the objective function and constraints by specifying the functions $\alpha$ and $\beta$.
\begin{itemize}
	\item Objective function: According to \eqref{eq:objective-general}, we can simply let $\alpha = f$ and $\beta=g$, as is shown in Row 2 of Table \ref{tab:unified-form}.
	\item Constraints: According to \eqref{eq:sc-constraint-tractable} and \eqref{eq:sc-constraint-var}, it is necessary to compute $\mathbb{E}^{\bm{x}_0,\bm{z}_0}\bm{\phi}_i^\top\bm{X}_t\bm{X}_t^\top\bm{\phi}_i$ and $\mathbb{E}^{\bm{x}_0,\bm{z}_0}\bm{\phi}_i^\top\bm{X}_t$. Rows 3 and 4 of Table \ref{tab:unified-form} show the specific SAF parameters for these 2 expectation values.
\end{itemize}

\begin{table}[!t]
	\renewcommand{\arraystretch}{1.1}
	\linespread{1.15}
	\centering
	\begin{small}
		\caption{SAF Form of the Objective Function and Constraints}
		\begin{tabular}{cccc}
			\hline\hline
			Type&$\mathcal{P}_{\alpha,\beta}(t,\bm{x}_0,\bm{z}_0)$&$\alpha$&$\beta$\\
			\hline
			Objective&$J$&$f$&$g$\\
			\hline
			Constraint&$\mathbb{E}^{\bm{x}_0,\bm{z}_0}\bm{\phi}_i^\top\bm{S}_t\bm{S}_t^\top\bm{\phi}_i$&0&$\bm{\phi}_i^\top\bm{S}_t\bm{S}_t^\top\bm{\phi}_i$\\
			\hline
			Constraint&$\mathbb{E}^{\bm{x}_0,\bm{z}_0}\bm{\phi}_i^\top\bm{S}_t$&0&$\bm{\phi}_i^\top\bm{S}_t$\\
			\hline \hline
		\end{tabular}\label{tab:unified-form}
	\end{small}
	\vspace{-3mm}
\end{table}

Based on the SAF formulation of the objective function and constraints, the SCP can be formulated as

\noindent\textbf{SCP}:
\begin{equation} \label{eq:scp}
\begin{split}
&\min_{u\in \mathcal{U}} J = \mathcal{P}_{f, g}(T, \bm{x}_0, \bm{z}_0)\\
&\text{s.t. } \\
&d\bm{X}_t = (\bm{A}\bm{X}_t + \bm{B}u(\bm{Z}_t) + \bm{C}\bm{Z}_t)dt\\
&d\bm{Z}_t = \mu(\bm{Z}_t) + \sigma(\bm{Z}_t)d\bm{W}_t \\
&\bm{X}_0=\bm{x}_0,\bm{Z}_0=\bm{z}_0\\
&\kappa_\gamma\sqrt{\mathcal{P}_{0,\bm{\phi}_i^\top\bm{S}_t\bm{S}_t^\top\bm{\phi}_i}(t, \bm{x}_0, \bm{z}_0) -\mathcal{P}_{0,\bm{\phi}_i^\top\bm{S}_t}(t,\bm{x}_0, \bm{z}_0)} \\
&\text{~~~}+ \mathcal{P}_{0,\bm{\phi}_i^\top\bm{S}_t}(t,\bm{x}_0, \bm{z}_0)\leq \bar{\phi}_i, \forall i \in \Omega^C
\end{split}
\end{equation}

Here, we highlight the differences between \eqref{eq:scp} and an MPC problem \cite{Ersdal2016Model} for an AGC system. The major difference is that the SCP explicitly considers the probability distribution of the uncertainty, which is characterized by $\sigma(\bm{Z}_t)d\bm{W}_t$. Because the SCP considers the stochastic characteristics of the wind power uncertainty, the controller $u$ obtained via the SCP is a closed-loop controller that can be applied over time. In contrast, the MPC problem models the uncertainty as a predicted input, and the result obtained via MPC can be used in only one time step.
%Therefore, the optimal control problem in MPC needs to be solved repeatedly during its control horizon, while the SCP need to be solved only once, after which the controller $u$ is applied without solving an optimal control problem.

To solve the SCP, a key challenge is to evaluate the SAFs in \eqref{eq:scp}. In conventional methods, stochastic terms are usually evaluated by means of scenario-based approaches. In these approaches \cite{Yang2016Stochastic}, the uncertainty is expressed as a set of scenarios, and the assessment problem for each scenario can be solved in a deterministic way. Once the value for each scenario has been obtained, the SAF is obtained by taking the average over all scenarios. However, scenario-based approaches are time-consuming, which limits their potential for use in online applications \cite{Chen2018Unified}. This is the motivation for this paper, i.e., to provide an effective way to evaluate the SAFs without scenario generation to allow the SCP to be solved efficiently.

\section{Convex Formulation of SCP-AGC } \label{section:convex}
\begin{figure}
	\centering
	\tikzstyle{format} = [draw, thin, fill=blue!20, minimum width=5em, minimum height=2em]
	\tikzstyle{format-red} = [draw, thin, fill=red!20, minimum width=5em, minimum height=2em]
	\tikzstyle{medium} = [ellipse, draw, thin, fill=green!20]
	
	\begin{tikzpicture}[node distance=2cm, auto,>=latex', thick]
	\path[->] node[medium] (scp) {SCP: Eq. \eqref{eq:scp}};
	\path[->] node[format, below left of=scp, xshift=-2em] (ss) {$\bm{S}_t$: Eq. \eqref{eq:ito-reduced}, \eqref{eq:s}}
	(scp) edge (ss);
	\path[->] node[format, below right of=scp, xshift=2em] (saf) {SAF: Eq. \eqref{eq:objective-general}, \eqref{eq:sc-constraint-var}, \eqref{eq:sap}} (scp) edge (saf);
	\path[->] node[format-red, below left of=ss, xshift=-0.5em] (s-tilde) {$\bm{\widetilde{S}}_t$: Eq. \eqref{eq:x-tilde}} (ss) edge (s-tilde);
	\path[->] node[format-red, below right of=ss, xshift=0.5em] (s-hat) {$\bm{\widehat{S}}_t$: Eq. \eqref{eq:x-hat}} (ss) edge (s-hat);
	\path[->] node[format-red, below of=saf, yshift=1.6em] (daf) {DAF: Eq. \eqref{eq:series-vn}} (saf) edge node {Theorem \ref{thm:series}}(daf);
	\path[->] node[medium, below of=scp, yshift=-7em] (scp-opt) {SCP-opt: Eq. \eqref{eq:scp-opt}} (s-tilde) edge (scp-opt) (s-hat) edge (scp-opt) (daf) edge (scp-opt);
%	(saf-pde) edge node {Theorem \ref{theorem:pde}} (daf-pde);

%	(daf-pde) edge node[near end] {Proposition \ref{prop:dap}}(daf);
	\end{tikzpicture}
	\caption{Flowchart of the proposed approach.}
	\label{fig:flowchart}
	\vspace{-5mm}
\end{figure}
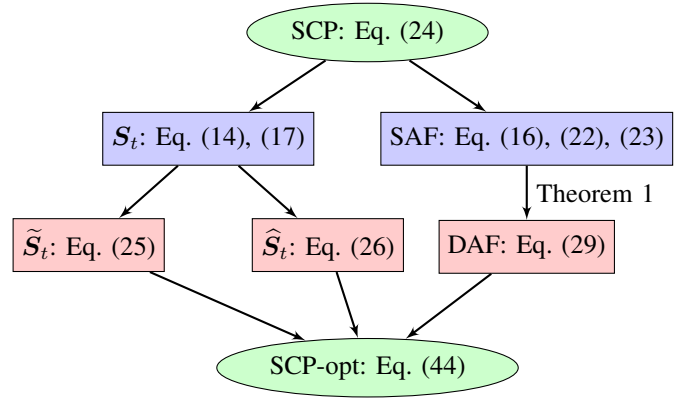

This section presents a deterministic convex optimization formulation of the SCP in \eqref{eq:scp}. The basic idea of this section is depicted in Fig. \ref{fig:flowchart}. A fundamental result presented in this section is the series expansion of an SAF, which makes it possible to transform the stochastic variables in $\bm{S}_t$ into 2 groups of deterministic variables, namely, $\bm{\widetilde{S}}_t$ and $\bm{\widehat{S}}_t$, and the SAFs into DAFs. $\bm{\widetilde{S}}_t$, $\bm{\widehat{S}}_t$, and the DAFs can be calculated in a deterministic way; thus, the SCP can be expressed as a deterministic optimization problem, denoted by SCP-opt. We then prove the convexity of SCP-opt and discuss the computational burden of the convex optimization problem.

\vspace{-2mm}
\subsection{Series Expansion of an SAF}
Before introducing the series expansion of an SAF, we first introduce two groups of auxiliary variables.

Let $\bm{\widetilde{X}}_t$ and $\bm{\widetilde{Z}}_t$ be defined as the solution to the following system of differential equations:
\vspace{-1mm}
\begin{equation} \label{eq:x-tilde}
\begin{split}
d\bm{\widetilde{X}}_t &= (\bm{A}\bm{\widetilde{X}}_t + \bm{B}u(\bm{\widetilde{Z}}_t) + \bm{C}\bm{\widetilde{Z}}_t)dt\\
d\bm{\widetilde{Z}}_t &= \mu(\bm{\widetilde{Z}}_t)dt\\
\bm{\widetilde{X}}_0&=\bm{x}_0,\bm{\widetilde{Z}}_0=\bm{z}_0\\
\end{split}
\end{equation}

\vspace{-1mm}
Let $\bm{\widehat{X}}_t$ and $\bm{\widehat{Z}}_t$ be defined as the solution to the following system of differential equations:
\vspace{-1mm}
\begin{equation} \label{eq:x-hat}
\begin{split}
d\bm{\widehat{X}}_t &= (\bm{A}\bm{\widehat{X}}_t + \bm{B}\frac{\partial u}{\partial \bm{\widetilde{Z}}_t}\bm{\widehat{Z}}_t + \bm{C}\bm{\widehat{Z}}_t)dt\\
d\bm{\widehat{Z}}_t &= \frac{\partial \mu}{\partial \bm{\widetilde{Z}}_t}\bm{\widehat{Z}}_tdt\\
\bm{\widehat{X}}_0&=\bm{0}_{N_x\times N_z},\bm{\widehat{Z}}_0=\bm{I}_{N_z}
\end{split}
\end{equation}
where $\bm{I}_{N_z}$ is the identity matrix of order $N_z$. Note that $\bm{\widehat{X}}_t$ is an $N_x\times N_z$ matrix, and $\bm{\widehat{Z}}_t$ is an $N_z \times N_z$ matrix.

For simplicity, let 
\vspace{-2mm}
\begin{equation} \label{eq:s-tilde-hat}
\bm{\widetilde{S}}_t = 
\begin{bmatrix}
\bm{\widetilde{X}}_t \\
u(\bm{\widetilde{Z}}_t) \\
\bm{\widetilde{Z}}_t
\end{bmatrix},
\bm{\widehat{S}}_t = 
\begin{bmatrix}
\bm{\widehat{X}}_t \\
\frac{\partial u}{\partial \bm{\widetilde{Z}}_t}\bm{\widehat{Z}}_t \\
\bm{\widehat{Z}}_t
\end{bmatrix}
\vspace{-1mm}
\end{equation}

It is clear that $\bm{\widehat{S}}_t = \partial \bm{\widetilde{S}}_t / \partial \bm{z}_0$.

%The relationship between $\bm{\widetilde{S}}_t$ and $\bm{\widehat{S}}_t$ is
%\begin{lemma} \label{lemma:s-tilde-hat}
%	$\bm{\widehat{S}}_t = \partial \bm{\widetilde{S}}_t / \partial \bm{z}_0$.
%\end{lemma}
%\textit{Proof}: Take the derivatives over $\bm{z}_0$ in each equation of \eqref{eq:x-tilde}, then we obtain the differential equation of $\partial \bm{\widetilde{S}}_t / \partial \bm{z}_0$, which is the same as \eqref{eq:x-hat}.
%\hfill $\square$

Notably, $\bm{\widetilde{S}}_t$ and $\bm{\widehat{S}}_t$ are both deterministic. In contrast, $\bm{S}_t$ is stochastic because $\bm{Z}_t$ is determined by an SDE, specifically the term $\sigma(\bm{Z}_t)d\bm{W}_t$. In practice, ordinary differential equations can be solved simply via discretization, while solving SDEs is much more difficult. The following theorem makes it possible to replace $\bm{S}_t$ with $\bm{\widetilde{S}}_t$ and $\bm{\widehat{S}}_t$ in the SCP.

\begin{theorem} \label{thm:series}
	Let $v(t,\bm{x}_0,\bm{z}_0)=\mathcal{P}_{\alpha,\beta}(t,\bm{x}_0,\bm{z}_0)$. If $\alpha$ and $\beta$ have continuous 2nd-order derivatives, then the following series is convergent
	\vspace{-1mm}
	\begin{equation} \label{eq:series}
	v(t,\bm{x}_0,\bm{z}_0)=\sum_{n=0}^\infty\widetilde{v}_n(t,\bm{x}_0,\bm{z}_0)
	\vspace{-1mm}
	\end{equation}
	with
	\begin{align}
	\widetilde{v}_n &= \int_0^t\alpha_n(\bm{\widetilde{S}}_s)ds + \beta_n(\bm{\widetilde{S}}_t) \label{eq:series-vn} \\
%	\begin{aligned}
	\alpha_0 &= \alpha, \beta_0=\beta \label{eq:series-0}\\
	\alpha_{n+1}(\bm{\widetilde{S}}_s)&=\frac{1}{2}\sigma(\bm{\widetilde{Z}}_s)^{\top}\left[\nabla^2_{\bm{z}_0}\widetilde{v}_n\right]\sigma(\bm{\widetilde{Z}}_s), n\geq 0 \label{eq:series-alpha-n}\\
	\beta_{n+1} &= 0, n \geq 0 \label{eq:series-beta-n} \\
	\left[\nabla^2_{\bm{z}_0}\widetilde{v}_n\right] &=\int_0^t \bm{\widehat{S}}_s^{\top}\left[\nabla^2\alpha_n\right]\bm{\widehat{S}}_sds +\bm{\widehat{S}}_t^{\top}\left[\nabla^2\beta_n\right]\bm{\widehat{S}}_t \label{eq:series-nabla-n}
	\end{align}
	where $\left[\nabla^2_{\bm{z}_0}\widetilde{v}_n\right]$ is the 2nd-order partial derivative of $\widetilde{v}_n$ with respect to $\bm{z}_0$, while $\nabla^2 \alpha_n$ and $\nabla^2 \beta_n$ are the Hessian matrices of $\alpha_n$ and $\beta_n$, respectively. The $\widetilde{v}_n$ are called DAFs in this paper, as shown in Fig. \ref{fig:flowchart}.
\end{theorem}

The proof of Theorem \ref{thm:series} is rather technical and is not closely related to the other content of this paper.
\begin{extended}
	Thus, we delay its presentation to Appendix \ref{appendix:proof-series}.
\end{extended}
\begin{reduced}
	We present its presentation in Appendix B of the extended version \cite{Chen2018Stochastic}.
\end{reduced}

According to Theorem \ref{thm:series}, an SAF can be expressed as the sum of a series of $\widetilde{v}_n$, $n \geq 0$. The stochastic term $\sigma(\cdot)$ does not exist in an SDE as in \eqref{eq:ito-reduced} but rather in a deterministic integral, as shown in \eqref{eq:series-alpha-n}.
%%%Editor - Please ensure that the intended meaning has been maintained
%%%in the above edit.
Moreover, since $\bm{\widetilde{S}}_t$ can be obtained in a deterministic way, each term in this series is a deterministic function of $t$, $\bm{x}_0$ and $\bm{z}_0$; i.e., no scenarios or expectation values need to be considered when computing this function. Therefore, Theorem \ref{thm:series} provides a means of computing the SAF in a deterministic way.

%In \eqref{eq:series-vn} $\sim$ \eqref{eq:series-beta-n}, $\bm{\widetilde{S}}_t$ can be obtained via \eqref{eq:x-tilde}. The only unsolved problem when computing $\mathcal{P}_{\alpha, \beta}$ is the computation of $\left[\nabla^2_{\bm{z}}\widetilde{v}_n\right]$ in \eqref{eq:series-alpha-n}. However, we have the following theorem:
%%%Editor - Please note that the sentence above appears to refer to a theorem
%%%that is commented out below and thus does not appear in the current version
%%%of the manuscript. Please consider deleting the sentence above.

%\begin{theorem} \label{thm:2nd}
%	If $\mu(\bm{Z}_t)$ is affine, then, the Hessian matrix $\left[\nabla^2_{\bm{z}_0}\widetilde{v}_n\right]$ in \eqref{eq:series-alpha-n} satisfies the following equality:
%	\begin{equation}
%	\left[\nabla^2_{\bm{z}_0}\widetilde{v}_n\right] =\int_0^t \bm{\widehat{S}}_s^{\top}\left[\nabla^2\alpha_n\right]\bm{\widehat{S}}_sds +\bm{\widehat{S}}_t^{\top}\left[\nabla^2\beta_n\right]\bm{\widehat{S}}_t
%	\end{equation}
%\end{theorem}
%where $\nabla^2 \alpha_n$ and $\nabla^2 \beta_n$ are the Hessian matrix of $\alpha_n$ and $\beta_n$ respectively.

By Theorem \ref{thm:series}, the computation of an SAF can be transformed into the computation of a series of DAFs, the computation of which depends only on $\bm{\widetilde{S}}_t$ and $\bm{\widehat{S}}_t$. Therefore, the equations for $\bm{S}_t$ can be removed from the SCP since they can be replaced with those for $\bm{\widetilde{S}}_t$ and $\bm{\widehat{S}}_t$.
%%%Editor - Please ensure that the intended meaning has been maintained
%%%in the above edit.

\vspace{-2mm}
\subsection{Reformulation of the Objective Function and Constraints} \label{section:scp-reformulation}
This subsection reformulates the objective function and constraints in accordance with Theorem \ref{thm:series}. Since the series expansion of an SAF is convergent, it can be approximated by a finite series. 
%For example, a zeroth-order approximation $v\approx \widetilde{v}_0$ of an SAF can be obtained by simply eliminating the diffusion, whereas a 1st-order approximation $v\approx \widetilde{v}_0+\widetilde{v}_1$ adds a 1st-order correction to the zeroth-order approximation.
Theoretically, this approximation can reach an arbitrary level of accuracy as long as sufficient high-order corrections are introduced. For convenience, we use the 1st-order approximation here as an example, but the same technique can be conveniently applied for higher orders of approximation.

\subsubsection{Objective Function}
According to Theorem \ref{thm:series}, we have
\begin{equation} \label{eq:d-obj}
J=\widetilde{J}_0+\widetilde{J}_1
\end{equation}
\begin{equation} \label{eq:d-j0}
\begin{split}
\widetilde{J}_0 = &\int_0^Tf(\bm{\widetilde{S}}_t)dt + g(\bm{\widetilde{S}}_T)
\end{split}
\end{equation}
\begin{align} \label{eq:d-j1}
\widetilde{J}_1 = \int_0^T\frac{1}{2}\sigma(\bm{\widetilde{Z}}_{t})^{\top}\bm{l}_{t}\sigma(\bm{\widetilde{Z}}_{t})dt
\end{align}
where $\bm{l}_t=\nabla^2_{\bm{z}_0}\widetilde{J}_0$.
%According to Theorem \ref{thm:series},
%%%Editor - Please note that the theorem referenced above does not appear in
%%%the manuscript.
$\bm{l}_{t}$ is calculated as follows:
\begin{align} \label{eq:d-l}
\bm{l}_{t} = \int_0^t\bm{\widehat{S}}_{s}^{\top}\left[\nabla^2f\right]\bm{\widehat{S}}_{s}ds+\bm{\widehat{S}}_{t}^{\top}\left[\nabla^2g\right]\bm{\widehat{S}}_{t}
\end{align}

To formulate a convex optimization problem, as discussed in Section \ref{section:optimization-convex}, we rewrite \eqref{eq:d-j0} and \eqref{eq:d-l} as inequalities:
\begin{equation} \label{eq:d-j0-geq}
\begin{aligned}
\widetilde{J}_0 \geq \int_0^Tf(\bm{\widetilde{S}}_t)dt + g(\bm{\widetilde{S}}_T)
\end{aligned}
\end{equation}
\begin{equation} \label{eq:d-l-geq}
\bm{l}_{t} \geq \int_0^t\bm{\widehat{S}}_{s}^{\top}\left[\nabla^2f\right]\bm{\widehat{S}}_{s}ds+\bm{\widehat{S}}_{t}^{\top}\left[\nabla^2g\right]\bm{\widehat{S}}_{t}
\end{equation}
It is clear that \eqref{eq:d-j0-geq} and \eqref{eq:d-l-geq} are equivalent to \eqref{eq:d-j0} and \eqref{eq:d-l} when $J$ is being minimized.

\subsubsection{Constraints}
For $\mathbb{E}^{\bm{x}_0,\bm{z}_0}\bm{\phi}_i^{\top}\bm{S}_{t}$, we have
\vspace{-1mm}
\begin{equation} \label{eq:constraint-1st}
\mathbb{E}^{\bm{x}_0,\bm{z}_0}\bm{\phi}_i^{\top}\bm{S}_{t}=\bm{\phi}_i^{\top}\mathbb{E}^{\bm{x}_0,\bm{z}_0}\bm{S}_{t}=\bm{\phi}_i^{\top}\bm{\widetilde{S}}_{t}
\vspace{-1mm}
\end{equation}

Now, let us consider $\mathbb{E}^{\bm{x}_0,\bm{z}_0}\bm{\phi}_i^{\top}\bm{S}_{t}\bm{S}_{t}^{\top}\bm{\phi}_i$. According to Theorem \ref{thm:series}, we have
\begin{equation}
\begin{split}
\mathbb{E}^{\bm{x}_0,\bm{z}_0}\bm{\phi}_i^{\top}\bm{S}_{t}\bm{S}_{t}^{\top}\bm{\phi}_i=&\bm{\phi}_i^{\top}\bm{\widetilde{S}}_{t}\bm{\widetilde{S}}_{t}^{\top}\bm{\phi}_i  \\ +&\int_0^t\sigma(\bm{\widetilde{Z}}_{s})^{\top}\bm{\widehat{S}}_{s}^{\top}\bm{\phi}_i\bm{\phi}_i^{\top}\bm{\widehat{S}}_{s}\sigma(\bm{\widetilde{Z}}_{s})ds
\end{split}
\end{equation}
in which we use the fact that $\nabla^2\left( \bm{\phi}_i^{\top}\bm{S}_{s}\bm{S}_{s}^{\top}\bm{\phi}_i\right)=\bm{\phi}_i\bm{\phi}_i^{\top}$.

According to \eqref{eq:sc-constraint-var}, we have
\vspace{-1mm}
\begin{equation}
\text{var}\left\{\bm{\phi}_i^{\top}\bm{S}_{t}\right\}=\int_0^t\sigma(\bm{\widetilde{Z}}_{s})^{\top}\bm{\widehat{S}}_{s}^{\top}\bm{\phi}_i\bm{\phi}_i^{\top}\bm{\widehat{S}}_{s}\sigma(\bm{\widetilde{Z}}_{s})ds
\vspace{-1mm}
\end{equation}

Thus, the constraint expressed in \eqref{eq:sc-constraint-tractable} can be rewritten as
\vspace{-1mm}
\begin{equation} \label{eq:d-constraint}
\kappa_\gamma\sqrt{\int_0^t\sigma(\bm{\widetilde{Z}}_{s})^{\top}\bm{\widehat{S}}_{s}^{\top}\bm{\phi}_i\bm{\phi}_i^{\top}\bm{\widehat{S}}_{s}\sigma(\bm{\widetilde{Z}}_{s})ds}+\bm{\phi}_i^{\top}\bm{\widetilde{S}}_{t}\leq \bar{\phi}_i
\vspace{-1mm}
\end{equation}

\vspace{-2mm}
\subsection{The Optimization Problem and Its Convexity}
\label{section:optimization-convex}
According to the above discussion, the optimization formulation of the SCP can be summarized as

\noindent\textbf{SCP-opt:}
\vspace{-2mm}
\begin{equation} \label{eq:scp-opt}
\begin{split}
\min_{u \in \mathcal{U}} J &:\eqref{eq:d-obj}\\
\text{s.t.}&: \eqref{eq:x-tilde}\eqref{eq:x-hat}\eqref{eq:d-j1}\eqref{eq:d-j0-geq}\eqref{eq:d-l-geq}\eqref{eq:d-constraint}
\end{split}
\end{equation}

Unlike in \eqref{eq:scp}, there is no stochastic term in SCP-opt. By discretizing SCP-opt over the time $t$, we can solve it with commercial optimization tools.

Moreover, regarding the convexity of the SCP-opt,
%, we have the following theorem:
%\begin{theorem} \label{t%hm:convex}
	if $f$ and $g$ are convex and of up to quadratic order and the control policy space $\mathcal{U}$ is convex, then SCP-opt is convex. %\footnote{Please do not mistake the convexity of the space $\mathcal{U}$ by the convexity of the function $u(\bm{Z}_t)$. The latter can be non-convex.}
%\end{theorem}
Actually, according to \eqref{eq:x-tilde} and \eqref{eq:x-hat}, $\bm{\widetilde{Z}}_{t}$ and $\bm{\widehat{Z}}_{t}$ are constants (i.e., independent of the function $u$). Therefore, \eqref{eq:x-tilde} and \eqref{eq:x-hat} are convex sets of $\bm{\widetilde{S}}_t$ and $\bm{\widehat{S}}_t$.
Since $f$ and $g$ are convex, \eqref{eq:d-j0-geq} is convex. Moreover, since $f$ and $g$ are convex and of up to quadratic order, $\nabla^2f$ and $\nabla^2g$ are constant positive semidefinite matrices; thus, the right-hand term of \eqref{eq:d-l-geq} is convex and quadratic, which means that \eqref{eq:d-l-geq} is convex.
Finally, \eqref{eq:d-constraint} is a second-order cone of $\bm{\widehat{S}}_{t}$ and $\bm{\widetilde{S}}_{t}$ and thus is convex. Therefore, SCP-opt is convex.

\begin{remark}
	Note that the convexity of $\mathcal{U}$ does not require the convexity of each $u \in \mathcal{U}$. Therefore, there are many ways to parametrize $\mathcal{U}$. For example, the affine disturbance control parametrization reads
	\vspace{-1mm}
	\begin{equation} \label{eq:u-affine}
	u(\bm{Z}_t) = \bm{U}_t^0 + \bm{F}_1^{\top}\bm{Z}_t
	\vspace{-1mm}
	\end{equation}
	while the quadratic disturbance control parametrization reads
	\vspace{-1mm}
	\begin{equation} \label{eq:u-quadratic}
	u(\bm{Z}_t) = \bm{U}_t^0 + \bm{F}_1^{\top}\bm{Z}_t + \bm{Z}_t^{\top}\bm{F}_2\bm{Z}_t
	\vspace{-1mm}
	\end{equation}
	In \eqref{eq:u-affine}, $\mathcal{U}$ is parametrized by $\bm{U}_t^0$ and $\bm{F}_1$, while in \eqref{eq:u-quadratic}, $\mathcal{U}$ is parametrized by $\bm{U}_t^0$, $\bm{F}_1$ and $\bm{F}_2$. Therefore, $\mathcal{U}$ is convex with regard to its parameters in each case, though in the 2nd case, $u\in\mathcal{U}$ is quadratic.
\end{remark}
\vspace{-2mm}
\subsection{Discussion of the Computational Burden}
SCP-opt is a convex optimization problem, which can be efficiently solved. Here, we discuss the computational burden of the proposed approach by showing that it is as efficient as the corresponding deterministic control (DC) problems and much faster than the MPC and SBSP approaches.
Note that the stochastic control problem degenerates into a deterministic control problem if $\sigma=0$,
and
therefore, the constraints \eqref{eq:d-j1} and \eqref{eq:d-l-geq} and the square root in \eqref{eq:d-constraint} do not exist. Since $\bm{\widehat{S}}_t$ appears only in these expressions, \eqref{eq:x-hat} is not needed. 
However, constraints \eqref{eq:x-tilde}, \eqref{eq:d-j0-geq} and \eqref{eq:d-constraint} are still needed to obtain a solution to the deterministic discrete control problem.
The numbers of variables and constraints are listed in Table \ref{tab:number-constraints},
where $N_x$ is the dimensionality of $X$, $N_z$ is the dimensionality of $\bm{Z}$, $N_c$ is the number of constraints, $N_t$ is the number of time steps in discretization, $N_p$ is the number of prediction steps in MPC, and $N_s$ is the number of scenarios in SBSP.

It is shown that the deterministic control problem contains $(N_x+N_z)N_t$ variables and $(N_x+N_z+N_c)N_t+1$ constraints,
whereas the proposed stochastic control approach involves $\left[(N_x+N_z)(N_z+1)+1\right]N_t$ variables and $\left[(N_x+N_z)(N_z+1)+N_c+1\right]N_t+2$ constraints.
In practice, we usually have $N_z\ll N_x$, i.e., the dominant sources of uncertainty, e.g., buses that are connected to wind generators, are much fewer than the number of state variables. Therefore, the numbers of variables and constraints in the considered problem are close to those in the DC problem.
By contrast, the computational burden of MPC in one control step is equal to that of the DC approach, where $N_t$ is replaced with the number of prediction steps $N_p$, and the total computational burden of MPC is equal to the one-step computational burden multiplied by the number of control steps, i.e., $N_t$. The SBSP approach requires approximately $N_s(N_x+N_z)N_t$ variables and $N_s[(N_x+N_z+N_c)N_t+1]$ constraints, as shown in Table \ref{tab:number-constraints}. 
Therefore, the scale of the MPC (or SBSP) problem is approximately $N_p$ (or $N_s$) times larger than that of the DC problem. Moreover, since the computational burden grows superlinearly with the numbers of variables and constraints,
the computational burden of the MPC and SBSP approaches will increase rapidly as $N_p$ and $N_s$ increase. In practice, $N_p$ and $N_s$ are usually large; thus, the proposed approach is much more efficient than the MPC and SBSP approaches.

In summary, the proposed approach allows non-Gaussian stochastic characteristics to be explicitly considered while incurring a computational burden similar to that of DC.

%\begin{table}[!t]
%	%	\vspace{-5mm}
%	\renewcommand{\arraystretch}{1.1}
%	%	\linespread{1.5}
%	\centering
%	\begin{small}
%		\caption{No. of Variables}
%%		\hspace{-5mm}
%		\begin{tabular}{cccc}
%			\hline\hline
%			Variable&ITB&DC&SBSP\\
%			\hline
%			$\bm{\widetilde{S}}_{t}$&$(N_x+N_z)N_t$&$(N_x+N_z)N_pN_t$&$N_s(N_x+N_z)N_t$\\
%			\hline
%			$\bm{\widehat{S}}_t$&$(N_x+N_z)N_zN_t$&0&0\\
%			\hline
%			$\bm{l}_{t}$&$N_t$&0&0\\
%			\hline \hline
%		\end{tabular}\label{tab:number-variables}
%	\end{small}
%	\vspace{-4mm}
%\end{table}

%\begin{table}[!t]
%	%	\vspace{-5mm}
%	\renewcommand{\arraystretch}{1.1}
%	%	\linespread{1.5}
%	\centering
%	\begin{small}
%		\caption{No. of Constraints}
%		\hspace{-5mm}
%		\begin{tabular}{cccc}
%			\hline\hline
%			Eq.&ITB&MPC&SBSP\\
%			\hline
%			\eqref{eq:x-tilde}&$(N_x+N_z)N_t$&$(N_x+N_z)N_pN_t$&$N_s(N_x+N_z)N_t$\\
%			\hline
%			\eqref{eq:x-hat}&$(N_x+N_z)N_zN_t$&0&0\\
%			\hline
%			\eqref{eq:d-j1}&$1$&0&0\\
%			\hline
%			\eqref{eq:d-constraint}&$N_cN_t$&$N_cN_pN_t$&$N_sN_cN_t$\\
%			\hline
%			\eqref{eq:d-j0-geq}&$1$&1&$N_s$\\
%			\hline
%			\eqref{eq:d-l-geq}&$N_t$&0&0\\
%			\hline \hline
%		\end{tabular}\label{tab:number-constraints}
%	\end{small}
%	\vspace{-4mm}
%\end{table}

\begin{table*}[!t]
	%	\vspace{-5mm}
	\renewcommand{\arraystretch}{1.1}
	%	\linespread{1.5}
	\centering
	\begin{small}
		\caption{No. of Variables and Constraints}
		\hspace{-5mm}
		\begin{tabular}{c|ccc|cccccc}
			\hline\hline
			&\multicolumn{3}{c|}{Variables}&\multicolumn{6}{c}{Constraints}\\
			\cline{2-10}
			&$\bm{\widetilde{S}}_{t}$&$\bm{\widehat{S}}_t$&$\bm{l}_{t}$&\eqref{eq:x-tilde}&\eqref{eq:x-hat}&\eqref{eq:d-j1}&\eqref{eq:d-constraint}&\eqref{eq:d-j0-geq}&\eqref{eq:d-l-geq}\\
			\hline
			SCP-opt&$(N_x+N_z)N_t$&$(N_x+N_z)N_zN_t$&$N_t$&$(N_x+N_z)N_t$&$(N_x+N_z)N_zN_t$&$1$&$N_cN_t$&$1$&$N_t$\\%&$N_s(N_x+N_z)N_t$\\
			\hline
			DC&$(N_x+N_z)N_t$&0&0&$(N_x+N_z)N_t$&0&0&$N_cN_t$&1&0\\
			\hline
			MPC&$(N_x+N_z)N_pN_t$&0&0&$(N_x+N_z)N_pN_t$&0&0&$N_cN_pN_t$&1&0\\
			\hline
			SBSP&$N_s(N_x+N_z)N_t$&0&0&$N_s(N_x+N_z)N_t$&0&0&$N_sN_cN_t$&$N_s$&0\\
			\hline \hline
		\end{tabular}\label{tab:number-constraints}
	\end{small}
	\vspace{-4mm}
\end{table*}

\vspace{-2mm}
\section{Case Study} \label{section:case}
This section presents a test case involving the IEEE 118-bus system \cite{IEEE118}. We first evaluate the effectiveness of the reformulation made in Section \ref{section:scp-reformulation} and then compare the proposed approach with conventional approaches.

\vspace{-2mm}
\subsection{Settings}
The case considered is the IEEE 118-bus system \cite{IEEE118}, which contains 54 generators and 186 lines. 
%The topology is shown in Fig. \ref{fig:case-B-topology}, 
The parameters of this system can be found in \cite{IEEE118}. The base power is 100~MVA. The parameters of the primary control loop are fixed, and the supplementary control output $P^{ref}_i$ is controllable. The frequency constraint is set to $50\pm0.1$ Hz. The time step is 1~s \cite{Chen2018Unified}, and the control horizon is 100~s. The stochastic sources are 6 wind farms on buses 6, 11, 18, 32, 55 and 69, and the nominal power of each generator is 100~MW. The power output of each wind generator is the sum of a predicted power output $P^{pred}_w$ and a fluctuating power output $\Delta P_w$. In accordance with \cite{Bludszuweit2008Statistical}, we assume that the wind power uncertainty satisfies the Laplace distribution $L(0,0.05)$, i.e.,
\vspace{-1mm}
\begin{equation} \label{eq:case-B-wind}
d\Delta P_w = -\Delta P_wdt + (0.1|\Delta P_w|+0.005)dW_t
\vspace{-1mm}
\end{equation}

To verify the effectiveness of the proposed method, the analytical results obtained using the proposed approach will be compared with the results obtained through a Monte Carlo simulation. In the Monte Carlo simulation, 1000 wind generation scenarios satisfying \eqref{eq:case-B-wind} were generated. 
%Fig. \ref{fig:case-B-wind} shows the predicted wind power and one of the actual wind power scenarios.

%\begin{figure}[!t]
%	\centering
%	\includegraphics[width=0.8\columnwidth]{ieee118.eps}
%	\caption{IEEE 118-bus system.}
%	\label{fig:case-B-topology}
%	\vspace{-3mm}
%\end{figure}

%\begin{figure}[!t]
%	\centering
%	\includegraphics[width=0.8\columnwidth]{wind_caseB.eps}
%	\caption{Wind power profile.}
%	\label{fig:case-B-wind}
%	\vspace{-3mm}
%\end{figure}
\vspace{-3mm}
\subsection{Accuracy of 1st-Order Reformulation}
\begin{figure}[!t]
	\centering
	\includegraphics[width=0.8\columnwidth]{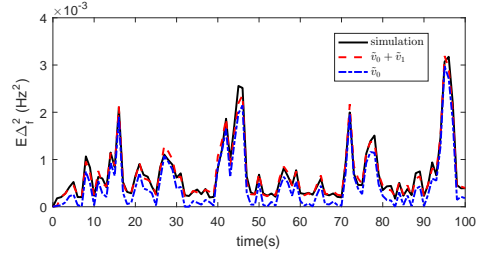}
	\vspace{-2mm}
	\caption{Simulated results for system frequency.}
	\label{fig:case-B-freq}
	\vspace{-3mm}
\end{figure}

In Section \ref{section:scp-reformulation}, we used a 1st-order reformulation of the SAF, i.e., $v\approx \widetilde{v}_0+\widetilde{v}_1$; in this subsection, we will validate the effectiveness of this order of reformulation for use in the SCP. A PI controller is used here, 
%i.e.
%\begin{equation}
%P^{ref}_i = -c_i\left(k_p\Delta_f+k_i\int \Delta_fdt\right),\forall i\in \Omega^G
%\end{equation}
%where $k_p$ and $k_i$ are the proportional and integral coefficients, and $c_i$ is the coefficient of power allocation satisfying $\sum_{i\in\Omega^G}c_i = 1$. Here we set $k_p=k_i=0.1$.
and the test SAF is the variance of the frequency deviation $\Delta_f(t)$ under wind power uncertainty. Therefore, $\alpha=0$ and $\beta=\Delta_f^2$. A simulation result obtained as an average value over 1000 Monte Carlo samples is also provided, and the variance of $\Delta_f$ is depicted in Fig. \ref{fig:case-B-freq}. It is evident that the first-order reformulation $\widetilde{v}_0+\widetilde{v}_1$ is in good agreement with the Monte Carlo simulation, while the zeroth-order reformulation $\widetilde{v}_0$ shows a larger deviation from the simulated value. This can be interpreted in terms of the concept of variance. In fact, the estimated value $\widetilde{v}_0$ exhibits frequency deviations based on the predicted wind power and is exactly $(\mathbb{E}\Delta_f)^2$, and the difference between $\mathbb{E}\Delta_f^2$ and $(\mathbb{E}\Delta_f)^2$ is the variance of $\Delta_f$, i.e., $\text{var}\left({\Delta_f}\right)$. Therefore, the difference between the first-order form and the zeroth-order form, i.e., $\widetilde{v}_1$, can be interpreted as the variance. 

In summary, the first-order reformulation achieves a good approximation of the SAF and captures the different characteristics of a stochastic system by means of different terms, i.e., $\widetilde{v}_0$ and $\widetilde{v}_1$.

%Moreover, the control performances of different control schemes can be assessed based on the proposed method. Fig. \ref{fig:case-B-control-scheme} shows the results for $\mathbb{E}\Delta_f^2$ obtained under different integral parameters $k_i$. A larger $k_i$ generally leads to a smaller frequency deviation, although exceptions arise at approximately 70s. These exceptions occur because the wind power continuously decreases, which results in a larger integral part of the control output. Therefore, a larger $k_i$ also leads to a larger overshoot.

%\begin{figure}[!t]
%	\centering
%	\includegraphics[width=0.7\columnwidth]{case-B-ki}
%	\caption{Frequency deviations under different control schemes.}
%	\vspace{-3mm}
%	\label{fig:case-B-control-scheme}
%\end{figure}
\vspace{-2mm}
\subsection{Performance and Benchmarks}
This section verifies the SCP-opt represented by \eqref{eq:scp-opt}. The objective function is given in \eqref{eq:objective}, where $\lambda_{ACE} = 10000, \bm{\Lambda}_U = 70000\bm{I}$, and $\mu_{ACE}=50000$. These weights are set to large values to ensure numerical convergence. Moreover, the affine control policy represented in \eqref{eq:u-affine} is used.
An initial control series $\bm{U}_t^0$ and a coefficient matrix $\bm{F}_1$ are computed and are then used to determine the control output in accordance with the actual wind power.

Here, we use the DC, MPC and PI approaches as benchmarks for conventional approaches. In the DC approach, an initial control series is computed that is not subsequently adjusted over the horizon. By contrast, in the MPC approach, the control policy is adjusted at each step by re-solving a finite-horizon deterministic control problem in accordance with the actual wind power. The MPC was designed to generate 10 steps (i.e., 10~s) of predicted values. In the PI controller, the parameters are tuned via a grid search. To illustrate the advantages of the proposed approach over the SBSP approach, SBSP schemes with 10 scenarios and 100 scenarios, denoted by SBSP(10) and SBSP(100), are also considered as benchmarks. In the SBSP approach, the control policy is the same as that in the SCP-opt approach, but the optimal control problem is solved via the scenario-based approach.

One thousand scenarios were generated to verify the performance of the proposed approach. The control performance is summarized in Table \ref{tab:case-B-comparison}, which shows that the SCP-opt and SBSP(100) achieve the best control performance. In comparison to the SCP-opt, the PI and DC approaches yield worse control performance in terms of both the objective function value and the probability that the frequency constraint will be violated. Specifically, since this case features high wind penetration, the PI controller cannot avoid constraint violation. The MPC approach achieves a performance similar to that of the proposed approach, although it is slightly worse because the number of prediction steps cannot be too large. SBSP(10) does not perform as well as the SCP-opt, whereas SBSP(100) performs slightly better. In fact, the models used in the SCP-opt and SBSP approaches are similar, but the SCP-opt uses a series expansion to solve the SCP, while the SBSP approach uses a scenario-based approach to solve it.

\begin{figure}[!t]
	\centering
	\includegraphics[width=0.9\columnwidth]{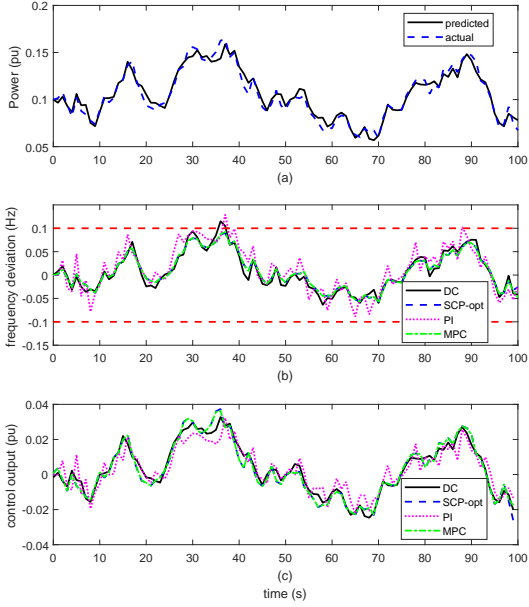}
	\vspace{-3mm}
	\caption{Comparison of 3 control approaches: (a) wind curves, (b) frequency deviations, and (c) control outputs.}
	\vspace{-3mm}
	\label{fig:case-B-comparison}
\end{figure}

Fig. \ref{fig:case-B-comparison} depicts the results of the SCP-opt, PI, DC and MPC approaches for one scenario. The results of the SBSP approach are not depicted because they are similar to those of the SCP-opt except for the computation time. Under the DC approach, the frequency exceeds the upper limit at 37 s because the actual wind power is higher than the predicted wind power. By contrast, in the SCP-opt and MPC approaches, the control output is adjusted to avoid this frequency violation. The control outputs of the SCP-opt and MPC approaches are similar, with the difference being that the control policy is recalculated just before 37 s in the MPC approach, whereas the control policy and the affine disturbance coefficient are calculated at the very beginning in the SCP-opt approach.

Now, we discuss the computation times, which are shown in Table \ref{tab:case-B-comparison}. The computation time of the PI controller is not reported because this controller does not perform an optimization process and thus is much faster than the other approaches. The computation time of the SCP-opt is approximately twice that of the DC approach. Although the average one-step computation time of the MPC approach is only 1.53~s because the number of prediction steps in the MPC scheme is smaller than the control horizons of the DC and SCP-opt, the total computation time of the MPC scheme is much greater than that of the SCP-opt because of its receding-horizon implementation. Moreover, the computation times for both SBSP(10) and SBSP(100) are unacceptable. This comparison clearly demonstrates the advantages of the proposed approach in terms of both performance and computational efficiency.

In summary, the proposed approach achieves good AGC performance in a computationally efficient manner. In contrast, the DC and PI approaches are computationally efficient but show worse performance, whereas the MPC and SBSP incur high computational burdens. Therefore, the proposed approach has attractive potential in online applications.

\begin{table}[!t] 
	%	\vspace{-5mm}
	\renewcommand{\arraystretch}{1.1}
	%	\linespread{1.5}
	\centering
	\begin{small}
		\caption{Simulation Results}
		%		\vspace{-2mm}
%		\begin{tabular}{c|ccccc}
%			\hline\hline
%			Approach&ITB&DC&MPC&SBSP(10)&SBSP(100)\\
%			\hline
%			Computation Time (s)&4.56&2.13&76.5&23.54&152.83\\
%			\hline
%			Objective Function&2838.4&2968.7&2884.7&2903.2&2832.5\\
%			\hline
%			Constraint Violation&0.7\%&29.7\%&0.5\%&1.2\%&0.7\%\\
%			\hline \hline
%		\end{tabular}
		\begin{tabular}{c|ccc}
			\hline\hline
			\multirow{2}{*}{Approach}&Computation&Objective &Constraint\\
			&Time (s)&Value&Violation (\%)\\
			\hline
			SCP-opt&8.62&2838.4&0.7\\
			\hline
			DC&4.26&2968.7&29.7\\
			\hline
			PI&-----&3571.3&100.0\\
			\hline
			MPC&153.1&2884.7&0.5\\
			\hline
			SBSP(10)&215.7&2893.2&1.2\\
			\hline
			SBSP(100)&1042.2&2832.5&0.7\\
			\hline \hline
		\end{tabular}
		\label{tab:case-B-comparison}
	\end{small}
	%	\vspace{-4mm}
\end{table}

\color{black}
\vspace{-3mm}
\section{Conclusion}\label{section:conclusion}
    This paper presents a systematic framework for the efficient optimal control of AGC systems. The It\^{o} model proposed in this paper is able to describe AGC systems with various non-Gaussian wind power uncertainty. Based on the It\^{o} model, a theorem of convergent series expansion of SAFs is proven, which enables the design of a highly efficient optimal stochastic control algorithm of It\^{o} AGC systems, and allows SCPs to be solved with a computational burden comparable to that for deterministic control problems without sacrificing performance.
    Computational applications show that the proposed approach outperforms some popular control algorithms while incurring affordable computational burden, and thus has attractive potential in the online applications of AGC systems.
    \vspace{-2mm}

\appendices
\begin{extended}

\section{The It\^{o} Coefficients of an Arbitrary Distribution} \label{appendix:ito-distribution}
Suppose that the PDF of $Z_t$ is $p(z)$. Then, we have the following proposition:
\begin{proposition} \label{prop:arbitrary-distribution}
	Suppose that the functions $\mu$ and $\sigma$ satisfy
	%	\vspace{-1mm}
	\begin{equation} \label{eq:relationship-mu-sigma}
	\begin{split}
	%	\mu(Z_t) &= -k(Z_t-r) \\
	\sigma^2(z) &= 2 \frac{\int_{-\infty}^{z}\mu(z')p(z')dz'}{p(z)}
	\end{split}
	\end{equation}
	where $p(\cdot)$ is a given probability density function (PDF).
%%%Editor - Please consider either defining PDF at its first appearance in the
%%%appendices (above) or omitting the definition here since it is already
%%%defined in the main text.
Then, the stationary PDF of $Z_t$ is $p$. 
\end{proposition}

\textit{Proof}: According to \eqref{eq:relationship-mu-sigma},
%\vspace{-1mm}
\begin{equation}
\frac{1}{2}\sigma^2(z)p(z) = \int_{\infty}^z\mu(z')p(z')dz'
\end{equation}
Taking the 2nd-order derivative with respect to $z$ on both sides yields
\begin{equation}
\frac{1}{2}\frac{\partial^2 \left(\sigma^2(z)p(z)\right)}{\partial z^2} = \frac{\partial \left(\mu(z)p(z)\right)}{\partial z}
\end{equation}
which is the steady-state Fokker-Planck equation \cite{Pardoux2014Stochastic} for the diffusion process of \eqref{eq:sc-source}. Therefore, $p(\cdot)$ is the stationary distribution of $Z_t$.
\hfill $\square$

\begin{corollary}
	The It\^{o} coefficients of the Gaussian distribution, the beta distribution, the gamma distribution, and the Laplace distribution are as listed in Table \ref{tab:ito-distribution}.
\end{corollary}
\section{Proof of Theorem \ref{thm:series}} \label{appendix:proof-series}
The proof of Theorem \ref{thm:series} is based on the partial differential equation (PDE) formulation of the SAF given in \eqref{eq:sap}.

\begin{lemma} \label{lemma:pde}
	Let $v(t,\bm{x}_0,\bm{z}_0)=\mathcal{P}_{\alpha,\beta}(t,\bm{x}_0,\bm{z}_0)$. Given a stochastic system $\bm{S}_t$ defined by \eqref{eq:ito-reduced} and \eqref{eq:s}, define the following auxiliary variables:
	\begin{equation*}
	b(\bm{x}_0, \bm{z}_0)=\begin{bmatrix}\bm{Ax}_0\\\bm{B}u(\bm{z}_0)+\bm{C}\bm{z}_0\end{bmatrix},\bm{s}_0=\begin{bmatrix}\bm{x}\\u(\bm{z})\\\bm{z}\end{bmatrix}
	\end{equation*}
	Then, the SAF $v$ satisfies the following equation:
	\begin{equation} \label{eq:pde-v}
	\begin{aligned}
	\frac{\partial v}{\partial t} =&
	b(\bm{x}_0,\bm{z}_0)^\top\frac{\partial v}{\partial \bm{x}_0}
	+ \mu(\bm{z}_0)^\top\frac{\partial v}{\partial \bm{z}_0} \\
	&+ \frac{1}{2}\sigma(\bm{z}_0)^\top\left[\nabla^2_{\bm{z}}v\right]\sigma(\bm{z}_0) + \alpha(\bm{s}_0)  \\
	v(0,\bm{x}_0, \bm{z}_0) =& \beta(\bm{x}_0, \bm{z}_0)
	\end{aligned}
	\end{equation}
\end{lemma}

This lemma is the famous Feynman-Kac formula. The proof is omitted here; instead, the reader is referred to \cite{Pardoux2014Stochastic}.

\begin{corollary} \label{corollary:vn}
	$\widetilde{v}_n$, as defined in \eqref{eq:series-vn}, satisfies the following PDE:
	\begin{equation} \label{eq:pde-vn}
	\begin{aligned}
	&
	\frac{\partial \widetilde{v}_n}{\partial t} =
	b(\bm{x},\bm{z})^\top\frac{\partial \widetilde{v}_n}{\partial \bm{x}}
	+ \mu(\bm{z})^\top\frac{\partial \widetilde{v}_n}{\partial \bm{z}} + \alpha_n(\bm{s})  \\
	&\widetilde{v}_n(0,\bm{x}, \bm{z}) = \beta_n(\bm{s})
	\end{aligned}
	\end{equation}
\end{corollary}
\noindent\textit{Proof}: According to \eqref{eq:series-vn}, the expression for $\widetilde{v}_n$ is similar to that for $\mathcal{P}_{\alpha,\beta}$ except that $\bm{S}_t$ is replaced by $\bm{\widetilde{S}}_t$. According to \eqref{eq:x-tilde} and \eqref{eq:s-tilde-hat}, the equations for $\bm{\widetilde{S}}_t$ differs from those for $\bm{S}_t$ only by the removal of the stochastic term $\sigma(\cdot)$;
%%%Editor - Please ensure that the intended meaning has been maintained
%%%in the above edit.
therefore, by removing $\sigma$ from \eqref{eq:pde-v}, we obtain \eqref{eq:pde-vn}.
\hfill $\square$

According to Lemma \ref{lemma:pde} and Corollary \ref{corollary:vn}, it suffices to prove that the $v$ described by \eqref{eq:pde-v} can be expressed as a sum of a series of the $\widetilde{v}_n$ described by \eqref{eq:pde-vn}. To achieve this goal, we present a parametrized PDE of $w(t,\bm{x}_0,\bm{z}_0;\epsilon)$:
\begin{equation} \label{eq:pde-w}
\begin{aligned}
\frac{\partial w}{\partial t} =&
b(\bm{x}_0,\bm{z}_0)^\top\frac{\partial w}{\partial \bm{x}_0}
+ \mu(\bm{z}_0)^\top\frac{\partial w}{\partial \bm{z}_0} \\
&+ \frac{1}{2}\epsilon\sigma(\bm{z}_0)^\top\left[\nabla^2_{\bm{z}}w\right]\sigma(\bm{z}_0) + \alpha(\bm{s}_0)  \\
w(0,\bm{x}_0, \bm{z}_0) =& \beta(\bm{x}_0, \bm{z}_0)
\end{aligned}
\end{equation}

Regarding $w$, we have Lemmas \ref{lemma:convergence} and \ref{lemma:w-vn}:
\begin{lemma} \label{lemma:convergence}
	Consider the Taylor expansion of $w$ in $\epsilon$:
	\begin{equation} \label{eq:w-taylor}
	w(t,\bm{x}_0,\bm{z}_0;\epsilon) = \sum_{n=0}^\infty\frac{1}{n!}\left[\left.\frac{\partial^nw}{\partial \epsilon^n} \right|_{\epsilon=0}\right]\epsilon^n
	\end{equation}
	then this Taylor expansion is convergent.
\end{lemma}
\noindent\textit{Proof}: it suffices to prove that the real function $w(t,\bm{x}_0,\bm{z}_0; \epsilon)$ has an analytical continuation on the complex plane. Supposing that $w$ and $\epsilon$ in \eqref{eq:pde-w} are complex numbers and letting $\epsilon = \epsilon^r+i\epsilon^i$ and $w=w^r+iw^i$, where $i$ is the imaginary unit, we have
\begin{equation} \label{eq:pde-complex}
\begin{aligned}
\frac{\partial w^r}{\partial t} =& b(\bm{x}_0,\bm{z}_0)^\top\frac{\partial w^r}{\partial \bm{x}_0} + \mu(\bm{z}_0)^\top\frac{\partial w^r}{\partial \bm{z}_0} \\
&+ \frac{1}{2}\epsilon^r\sigma(\bm{z}_0)^\top\left[\nabla^2_{\bm{z}_0}w^r\right]\sigma(\bm{z}_0) \\
&- \frac{1}{2}\epsilon^i\sigma(\bm{z}_0)^\top\left[\nabla^2_{\bm{z}_0}w^i\right]\sigma(\bm{z}_0)+ \alpha(\bm{s}_0) \\
\frac{\partial w^i}{\partial t} =& b(\bm{x}_0,\bm{z}_0)^\top\frac{\partial w^i}{\partial \bm{x}_0} + \mu(\bm{z}_0)^\top\frac{\partial w^i}{\partial \bm{z}_0} \\
&+ \frac{1}{2}\epsilon^r\sigma(\bm{z}_0)^\top\left[\nabla^2_{\bm{z}_0}w^i\right]\sigma(\bm{z}_0) \\
&+ \frac{1}{2}\epsilon^i\sigma(\bm{z}_0)^\top\left[\nabla^2_{\bm{z}_0}w^r\right]\sigma(\bm{z}_0) \\
w^r(0,&\bm{x}_0,\bm{z}_0;\epsilon) = \beta(\bm{x}_0,\bm{z}_0)\\
w^i(0,&\bm{x}_0,\bm{z}_0;\epsilon) = 0
\end{aligned}
\end{equation}

Let $\delta_1 = \partial w^r/\partial \epsilon^r - \partial w^i/\partial \epsilon^i$ and $\delta_2 = \partial w^i/\partial \epsilon^r+\partial w^r/\partial \epsilon^i$; then, by taking derivatives with respect to $\epsilon^r$ and $\epsilon^i$ in \eqref{eq:pde-complex}, we find that $\delta_1$ and $\delta_2$ satisfy the following PDEs:
\begin{equation} \label{eq:pde-delta}
\begin{aligned}
\frac{\partial \delta_1}{\partial t} =& b(\bm{x}_0,\bm{z}_0)^\top\frac{\partial \delta_1}{\partial \bm{x}_0} + \mu(\bm{z}_0)^\top\frac{\partial \delta_1}{\partial \bm{z}_0}\\
&+ \frac{1}{2}\epsilon^r\sigma(\bm{z}_0)^\top\left[\nabla^2_{\bm{z}_0}\delta_1\right]\sigma(\bm{z}_0) \\
&- \frac{1}{2}\epsilon^i\sigma(\bm{z}_0)^\top\left[\nabla^2_{\bm{z}_0}\delta_2\right]\sigma(\bm{z}_0)\\
\frac{\partial \delta_2}{\partial t} =& b(\bm{x}_0,\bm{z}_0)^\top\frac{\partial \delta_2}{\partial \bm{x}_0} + \mu(\bm{z}_0)^\top\frac{\partial \delta_2}{\partial \bm{z}_0} \\
&+ \frac{1}{2}\epsilon^r\sigma(\bm{z}_0)^\top\left[\nabla^2_{\bm{z}_0}\delta_2\right]\sigma(\bm{z}_0) \\
&+ \frac{1}{2}\epsilon^i\sigma(\bm{z}_0)^\top\left[\nabla^2_{\bm{z}_0}\delta_1\right]\sigma(\bm{z}_0) \\
\delta_1(0,&\bm{x}_0,\bm{z}_0;\epsilon) = 0\\
\delta_2(0,&\bm{x}_0,\bm{z}_0;\epsilon) = 0
\end{aligned}
\end{equation}

It is clear that \eqref{eq:pde-delta} is a 2nd-order PDE with an initial condition of zero. Therefore, $\delta_1=\delta_2=0$, which means that
\begin{equation}
\begin{split}
\frac{\partial w^r}{\partial \epsilon^r} = \frac{\partial w^i}{\partial \epsilon^i} \\
\frac{\partial w^i}{\partial \epsilon^r} = - \frac{\partial w^r}{\partial \epsilon^i}
\end{split}
\end{equation}

Therefore, according to the Cauchy-Riemann condition, $w$ is holomorphic over the whole complex plane; thus, the radius of convergence of its Taylor expansion is $\infty$.
\hfill $\square$

\begin{lemma} \label{lemma:w-vn}
	The following equality holds for $\widetilde{v}_n$ and $w$:
	\begin{equation}
	\widetilde{v}_n=\frac{1}{n!}\left.\frac{\partial^nw}{\partial \epsilon^n} \right|_{\epsilon=0}
	\end{equation}
\end{lemma}
\noindent\textit{Proof}: when $n=0$, this lemma is trivial. Now, we assume that the lemma holds for $\widetilde{v}_n$; it is necessary to prove that the lemma holds for $\widetilde{v}_{n+1}$.

Taking the $(n+1)$-th derivative with respect to $\epsilon$ in \eqref{eq:pde-w} yields
\begin{equation} \label{eq:proof-pde-nth}
\begin{split}
\frac{\partial}{\partial t} \frac{\partial^{n+1} w}{\partial \epsilon^{n+1}} =&
b(\bm{x}_0,\bm{z}_0)^\top\frac{\partial}{\partial \bm{x}_0} \frac{\partial^{n+1} w}{\partial \epsilon^{n+1}} +
\mu(\bm{z}_0)^\top\frac{\partial}{\partial \bm{z}_0} \frac{\partial^{n+1} w}{\partial \epsilon^{n+1}} \\ &+\frac{1}{2}\sigma(\bm{z}_0)^\top\left[\nabla^2_{\bm{z}_0}\frac{\partial^{n+1}}{\partial \epsilon^{n+1}}(\epsilon w)\right]\sigma(\bm{z}_0)
\end{split}
\end{equation}

Consider the following equality:
\begin{equation} \label{eq:derivative-equality}
\frac{\partial^{n+1}(\epsilon w)}{\partial \epsilon^{n+1}} = (n+1)\frac{\partial^n w}{\partial \epsilon^n} + \epsilon \frac{\partial^{n+1}w}{\partial \epsilon^{n+1}}
\end{equation}

By applying \eqref{eq:derivative-equality} to \eqref{eq:proof-pde-nth} and taking $\epsilon=0$, we obtain
\begin{equation} \label{eq:pde-vn-derive}
\begin{split}
&\frac{\partial}{\partial t} \left[\frac{1}{(n+1)!}\left.\frac{\partial^{n+1} w}{\partial \epsilon^{n+1}}\right|_{\epsilon=0}\right] \\
=&
b(\bm{x}_0,\bm{z}_0)^\top\frac{\partial}{\partial \bm{x}_0} \left[\frac{1}{(n+1)!}\left.\frac{\partial^{n+1} w}{\partial \epsilon^{n+1}}\right|_{\epsilon=0}\right]\\ &+\mu(\bm{z}_0)^\top\frac{\partial}{\partial \bm{z}_0} \left[\frac{1}{(n+1)!}\left. \frac{\partial^{n+1} w}{\partial \epsilon^{n+1}}\right|_{\epsilon=0}\right] \\ &+\frac{1}{2}\sigma(\bm{z}_0)^\top\nabla^2_{\bm{z}_0}\left[\frac{1}{n!}\left.\frac{\partial^{n} w}{\partial \epsilon^{n}}\right|_{\epsilon=0}\right]\sigma(\bm{z}_0)
\end{split}
\end{equation}

Moreover, by applying \eqref{eq:series-alpha-n} and \eqref{eq:series-beta-n} to \eqref{eq:pde-vn}, we obtain
\begin{equation} \label{eq:pde-vn-exact}
\begin{aligned}
\frac{\partial \widetilde{v}_{n+1}}{\partial t} =&
b(\bm{x}_0,\bm{z}_0)^\top\frac{\partial \widetilde{v}_{n+1}}{\partial \bm{x}_0}
+ \mu(\bm{z}_0)^\top\frac{\partial \widetilde{v}_{n+1}}{\partial \bm{z}_0} \\
&+ \frac{1}{2}\sigma(\bm{z}_0)^\top\left[\nabla^2_{\bm{z}_0}\widetilde{v}_n\right]\sigma(\bm{z}_0)
\end{aligned}
\end{equation}

By comparing \eqref{eq:pde-vn-exact} with \eqref{eq:pde-vn-derive} and recalling that Lemma \ref{lemma:w-vn} holds for $n$, we conclude that
\begin{equation}
\widetilde{v}_{n+1} = \frac{1}{(n+1)!}\left.\frac{\partial^{n+1}w}{\partial \epsilon^{n+1}}\right|_{\epsilon=0}
\end{equation}

Therefore, Lemma \ref{lemma:w-vn} also holds for $n+1$, which completes the proof.
\hfill $\square$

\noindent\textit{Proof of Theorem \ref{thm:series}}:
it is clear that $v(t,\bm{x}_0,\bm{z}_0)=w(t,\bm{x}_0,\bm{z}_0; 1)$; therefore, by setting $\epsilon=1$ in \eqref{eq:w-taylor}, we obtain
\begin{equation}
v(t,\bm{x}_0,\bm{z}_0) = \sum_{n=0}^\infty\frac{1}{n!}\left[\left.\frac{\partial^nw}{\partial \epsilon^n} \right|_{\epsilon=0}\right]
\end{equation}
Then, according to Lemma \ref{lemma:w-vn}, we have $v=\sum_{n=0}^\infty\widetilde{v}_n$, and according to Lemma \ref{lemma:convergence}, the series is convergent. 

Formula \eqref{eq:series-nabla-n} remains to be proven. According to \eqref{eq:series-vn},
\begin{equation}
\left[\nabla^2_{\bm{z}_0}\widetilde{v}_n\right] =\int_0^t\left[\nabla^2_{\bm{z}_0}\alpha_n(\bm{\widetilde{S}}_s)\right]ds + \left[\nabla^2_{\bm{z}_0}\beta_n(\bm{\widetilde{S}}_t)\right]
\end{equation}

When $\mu$ is affine, we have
\begin{equation}
\left[\nabla^2_{\bm{z}_0}\alpha_n(\bm{\widetilde{S}}_s)\right] = \left[\frac{\partial \bm{\widetilde{S}}_s}{\partial \bm{z}_0}\right]^\top \left[\nabla^2 \alpha_n\right]\left[\frac{\partial \bm{\widetilde{S}}_s}{\partial \bm{z}_0}\right]
%=\bm{\widehat{S}}_s^{\top}\left[\nabla^2\alpha_n\right]\bm{\widehat{S}}_s
\end{equation}

Therefore,
\begin{equation}
\left[\nabla^2_{\bm{z}_0}\alpha_n(\bm{\widetilde{S}}_s)\right] =\bm{\widehat{S}}_s^{\top}\left[\nabla^2\alpha_n\right]\bm{\widehat{S}}_s
\end{equation}

A similar procedure can be applied to $\nabla^2_{\bm{z}_0}\beta_n(\bm{\widetilde{S}}_t)$. Thus, we complete the proof of Theorem \ref{thm:series}.
\hfill$\square$
%\section{Proof of Proposition \ref{prop:convergence}} \label{appendix:convergence}

%\section{Proof of Proposition \ref{prop:w-vn}} \label{appendix:w-vn}

\end{extended}

\bibliographystyle{ieeetr}
\bibliography{agc-scp}

\end{document}